\documentclass{ifacconf}

\usepackage{graphicx}  
\usepackage{amsmath,amsfonts}  
\usepackage{amssymb}
\usepackage{mathtools}
\usepackage{natbib}  
\usepackage{tikz}
\usepackage[tikz]{bclogo}

\newcommand{\R}{\mathbb{R}}
\newcommand{\C}{\mathbb{C}}

\newcommand{\F}{\mathbb F}
\newcommand{\A}{\mathcal{A}}

\newcommand{\BB}{\mathcal{B}}
\newcommand{\CC}{\mathcal{C}}
\newcommand{\HH}{\mathcal{H}}
\newcommand{\G}{\mathcal{G}}
\renewcommand{\Re}{\mathrm{Re}}
\renewcommand{\Im}{\mathrm{Im}}

\begin{document}
\begin{frontmatter}

\title{Well-posedness of a class of infinite-dimensional port-Hamiltonian systems with boundary control and observation\thanksref{footnoteinfo}} 

\thanks[footnoteinfo]{Bouchra Elghazi acknowledges funding from the European Union Horizon Europe MSCA Grant No. 101073558 (ModConFlex).}

\author[1]{Bouchra Elghazi}
\author[1]{Birgit Jacob}  
\author[2,3]{Hans Zwart}

\address[1]{School of Mathematics and Natural Sciences, University of Wuppertal, Gaußstraße 20, 42119 Wuppertal, Germany, \\ (e-mail: elghazi@uni-wuppertal.de, bjacob@uni-wuppertal.de).}
\address[2]{Department of Applied Mathematics, University of Twente, P.O. Box 217, 7500 AE, Enschede, The Netherlands, \\ (e-mail: h.j.zwart@utwente.nl)}
\address[3]{Department of Mechanical Engineering, Eindhoven University of Technology, P.O. Box 513, 5600 MB Eindhoven,The Netherlands, (e-mail: h.j.zwart@tue.nl)}

\begin{abstract} 
We characterize the well-posedness of a class of infinite-dimensional port-Hamiltonian systems with boundary control and observation. This class includes in particular the Euler-Bernoulli beam equations and more generally 1D linear infinite-dimensional port-Hamiltonian systems with boundary control and observation as well as coupled systems. It is known, that for the Timoshenko beam models internal well-posedness implies well-posedness of the overall system. By means of an example we show that this is not true for the Euler-Bernoulli beam models. An easy verifiable equivalent condition for well-posedness of the overall system will be presented. We will conclude the paper by applying the obtained results to several Euler-Bernoulli beam models.
\end{abstract}


\end{frontmatter}

\section{Introduction}

Port-Hamiltonian systems (pHS) provide a unified energy-based perspective on a wide range of physical systems, allowing a structured approach to their analysis, design, and control (\cite{SchMas:02,Schaft:06,SchJelt:14}). In this paper, we consider linear distributed port-Hamiltonian systems on a one-dimensional spatial domain with boundary control and boundary observation 
\begin{equation*}
\begin{aligned}
    \frac{\partial x}{\partial t}(\xi,t) &= \sum_{k=0}^{N} P_k \frac{\partial^k}{\partial \xi^k} \HH(\xi) x(\xi,t) , \quad \xi \in [a,b], \\
    u(t) &= W_{B,1} \tau(\HH x) \ , \ 0 = W_{B,2} \tau(\HH x), \\ 
    y(t) &= W_C \tau(\HH x) \ ,
\end{aligned}
\end{equation*}
where $\tau : H^N([a,b];\F^n) \to \F^{2 N n}$ is the trace operator, given by
\begin{equation*}
    \tau(x) =
    \begin{pmatrix}
        x(b) & \cdots & x^{(N-1)}(b) & x(a) & \cdots & x^{(N-1)}(a)
    \end{pmatrix}^\top.
\end{equation*}
The well-posedness of such systems has mainly been studied in the specific case $N=1$. In this context, it has been shown that the internal well-posedness, which means that the main operator generates a $C_0$-semigroup, is equivalent to the well-posedness of the system (\cite{ZwGorMasVill:10,JacZw:12}). However, this result does not extend to the case where $N>1$, see Example \ref{example}. Whereas the Timoshenko beam models are covered by the case $N=1$, the Euler-Bernoulli beam models require $N=2$. In this article we investigate the well-posedness for $N=2$ under the assumption of $P_1=0$ and $\HH$ being constant. More precisely, we investigate port-Hamiltonian systems of the following form
\begin{equation}\label{eqn:system1}
\begin{aligned}
    \frac{\partial x}{\partial t}(\xi,t) &= \left( P_2 \frac{\partial^2}{\partial \xi^2} + P_0 \right)\HH x(\xi,t), \quad \xi \in [a,b], \\
     u(t) &= \!W_{B,1} \!
     \begin{bmatrix}
        \HH x(b,t) \\ \HH x'(b,t) \\
        \HH x(a,t) \\  \HH x'(a,t)
    \end{bmatrix},\,   
    0 = \!W_{B,2} \!
    \begin{bmatrix} 
        \HH x(b,t) \\ \HH x'(b,t) \\ 
        \HH x(a,t)\\  \HH x'(a,t)
    \end{bmatrix},\, \\
      y(t) &= \!W_C \!\begin{bmatrix}
        \HH x(b,t)\\  \HH x'(b,t)\\
        \HH x(a,t)\\  \HH x'(a,t)
    \end{bmatrix},\, \\
    x(\xi,0) &= x_0(\xi), \quad \xi \in [a,b],
\end{aligned}
\end{equation}
where $t\ge 0$ denotes the time, $x(\cdot,t)\in \mathbb F^n$ is the state, $u(t)\in \mathbb F^m$ the input, $y(t)\in \mathbb F^m$ the output at time $t$ and $x'$ denotes the derivative of $x$ w.r.t $\xi$. Here $\mathbb F\in \{\mathbb R, \mathbb C\}$.

We assume that $P_0$, $P_2$, $\HH \in \mathbb F^{n\times n}$, $P_0$ and $P_2$ are skew-adjoint, $P_2$ is invertible and $\HH$ is self-adjoint and positive, that is, $x^\ast \HH x >0$ for $x\not=0$. Further, $W_{B,1}, W_C \in \F^{m \times 4n}$ with $0<m\le 2n$, and $W_{B,2}\in \F^{(2n-m) \times 4n}$ such that $\left[\begin{smallmatrix} W_{B,1}\\ W_{B,2} \\ W_C \end{smallmatrix}\right]$ has full row rank and
\begin{equation}\label{eqn:dissi}
\Re \langle P_2 \HH x'',\HH x\rangle_{L^2} \le  2\Re (W_{B,1}\HH\tau(x))^\ast W_{C}\HH \tau(x) , 
\end{equation}
for $x\in H^2([a,b],\F^n)$ satisfying $W_{B,2}\HH \tau(x)=0$. Here 
$H^2([a,b],\F^n)$ denotes the second Sobolev space.

If $m=2n$, then \eqref{eqn:dissi} is equivalent to 
\begin{equation*}
    \begin{bmatrix}  \tilde W_B\Sigma \tilde W_B^\ast &  \tilde W_B\Sigma \tilde W_C^\ast \\  \tilde W_C\Sigma \tilde W_B^\ast &  \tilde W_C\Sigma \tilde W_C^\ast \end{bmatrix}^{-1}\le \begin{bmatrix}  0 & I \\ I & 0 \end{bmatrix},
\end{equation*}

 where $\Sigma=\left[\begin{smallmatrix}  0 & I \\ I & 0 \end{smallmatrix}\right]$, $\tilde W_B= \frac{1}{\sqrt{2}} \left[\begin{smallmatrix}
     W_{B,1} \\ W_{B,2}
 \end{smallmatrix}\right] \left[\begin{smallmatrix}
       R & I \\
      -R & I
 \end{smallmatrix} \right]$, $\tilde W_C= \frac{1}{\sqrt{2}} W_C \left[\begin{smallmatrix}
       R & I \\
      -R & I
 \end{smallmatrix} \right]$ and $R= \left[\begin{smallmatrix}
     0 & -P_2^{-1} \\
     P_2^{-1} & 0
 \end{smallmatrix} \right]$, see \cite{Aug:16}.

System \eqref{eqn:system1} is equipped with the Hamiltonian
\begin{equation*}
H(x(\cdot))= \frac{1}{2}\int_a^b x(\xi)^\ast \HH x(\xi) d\xi, 
\end{equation*}
where $x\in L^2([a,b],\mathbb F^n)$. The port-Hamiltonian system satisfies the \emph{dissipativity inequality}, that is,
   \begin{equation*}
        \frac{d}{dt} H (x(\cdot,t)) \leq  \Re\,  u(t)^\ast y(t) 
    \end{equation*}
for classical solution of system \eqref{eqn:system1}, see \cite{Aug:16}. 

Using the port-Hamiltonian approach, internal well-posed\-ness (\cite{GorZwaMas:2005,Villegas:07}), stability and stabilizability (\cite{AugJac:14,Aug:15,Aug:16}) and robust output regulation (\cite{HumLassiPohj:16,HumLassi:18}) of the system \eqref{eqn:system1}
have been investigated. In this paper, we characterize the well-posedness of the system \eqref{eqn:system1} with boundary control and boundary observation by means of an easy verifiable matrix condition.  

\section{Example class: Euler-Bernoulli beam equation}

The Euler-Bernoulli beam model of a beam of length 1 is given by
\begin{equation}\label{eqn:beam}
 \rho \frac{\partial^2 w}{\partial t^2}(\xi,t) + \frac{\partial^2}{\partial \xi^2}\left( EI\frac{\partial^2 w}{\partial \xi^2}(\xi,t)\right)=0, 
\end{equation} 
where $\xi \in (0,1)$ and $t \geq 0$. Here $w(\xi, t)$ denotes the transverse vibrations of the beam at position $\xi \in (0, 1)$ and time $t\ge 0$. Further, $\rho$ is the mass per unit length, $E$ the Young’s modulus of the beam, and $I$ the area moment of inertia of the beam’s cross section. 

 Next, we aim to write the partial-differential equation in the form 
 $$   \frac{\partial x}{\partial t}(\xi,t) = \left(P_2 \frac{\partial^2}{\partial \xi^2} +P_0\right)\HH x(\xi,t) $$
 by a suitable choice of the state variable $x$.
 Defining 
 $$x= 
    \begin{bmatrix}
        \rho \tfrac{\partial w}{\partial t} \\[0.8ex] 
        \tfrac{\partial^2 w}{\partial \xi^2} 
    \end{bmatrix},$$
the Euler-Bernoulli beam equations \eqref{eqn:beam} is equivalent to 
   \begin{equation*}
    \frac{\partial x}{\partial t} (\xi,t) = 
    \begin{bmatrix} 
        0 & -1 \\ 1 & 0 
    \end{bmatrix} 
    \frac{\partial^2}{\partial \xi^2}
    \begin{bmatrix}
        \tfrac{1}{\rho} & 0 \\ 0 & EI
    \end{bmatrix} x(\xi,t).
   \end{equation*} 
So we set
   $$ P_2 = \begin{bmatrix} 0 & -1 \\ 1 & 0 \end{bmatrix}, \quad P_0=0, \quad 
   \HH = \begin{bmatrix}
        \tfrac{1}{\rho} & 0 \\ 0 & EI
    \end{bmatrix}. 
    $$
If we define the boundary control and observation by
\begin{equation}\label{eqn:uy}
    \begin{bmatrix}
        0 \\ I
    \end{bmatrix}\! u(t) \!=\! \begin{bsmallmatrix}
       \! W_{B,1} \!\\\! W_{B,2}\!
    \end{bsmallmatrix}\!\!
    \begin{bsmallmatrix}
        \tfrac{\partial w}{\partial t}(1,t) \\[0.8ex]  
        \!EI \tfrac{\partial^2 w}{\partial \xi^2}(1,t)\! \\[0.8ex] 
        \tfrac{\partial^2 w}{\partial t\partial \xi} (1,t)\\[0.8ex]  
        \!EI \tfrac{\partial^3 w}{\partial \xi^3} (1,t)\!\\[0.8ex] 
        \tfrac{\partial w}{\partial t}(0,t) \\[0.8ex]  
        \!EI \tfrac{\partial^2 w}{\partial \xi^2}(0,t)\! \\[0.8ex] 
        \tfrac{\partial^2 w}{\partial t\partial \xi} (0,t)\\[0.8ex]  
        \!EI \tfrac{\partial^3 w}{\partial \xi^3} (0,t)\!
   \end{bsmallmatrix}\!,
   \,
   y(t) \!=\! W_C \!\begin{bsmallmatrix}
        \tfrac{\partial w}{\partial t}(1,t) \\[0.8ex]  
        \!EI \tfrac{\partial^2 w}{\partial \xi^2}(1,t) \!\\[0.8ex] 
        \tfrac{\partial^2 w}{\partial t\partial \xi} (1,t)\\[0.8ex]  
        \!EI \tfrac{\partial^3 w}{\partial \xi^3} (1,t)\!\\[0.8ex] 
        \tfrac{\partial w}{\partial t}(0,t) \\[0.8ex]  
        \!EI \tfrac{\partial^2 w}{\partial \xi^2}(0,t)\! \\[0.8ex] 
        \tfrac{\partial^2 w}{\partial t\partial \xi} (0,t)\\[0.8ex]  
        \!EI \tfrac{\partial^3 w}{\partial \xi^3} (0,t)\!
   \end{bsmallmatrix}  
\end{equation}  
and assume that the matrices $\begin{bsmallmatrix}
       W_{B,1} \\ W_{B,2}
   \end{bsmallmatrix}$, $W_C\in \F^{4\times 8}$ satisfy the assumptions made in the introduction, then boundary controlled and boundary observed Euler-Bernoulli beam equation fits into the class studied in this paper.

   Well-possedness of the boundary controlled and boundary observed Euler-Bernoulli beam equation means loosely speaking the following: for every initial condition $x(\cdot,0)\in L^2([0,1]; \F^4)$ and every control function $u\in L^2_{\mathrm{loc}}([0,\infty); \mathbb F^4)$ the solution of \eqref{eqn:beam} control and observation given by \eqref{eqn:uy} possess a mild solution $x$ and the output function satisfies $y\in L^2_{\mathrm{loc}}([0,\infty); \mathbb F^4)$. Further, the state and output depend continuously on the initial condition and input function. For the precise definition we refer the reader to the following section. We derive condition in terms of the matrices $\begin{bsmallmatrix}
       W_{B,1} \\ W_{B,2}
   \end{bsmallmatrix}$ and $W_C$ guaranteeing well-posedness of the Euler-Bernoulli beam.

\section{Well-posed boundary control systems}

We use an operator theoretical approach to study well-posedness of system \eqref{eqn:system1}. Therefore, we first rewrite system \eqref{eqn:system1} equivalently into an abstract boundary control and boundary observation system. 

\begin{defn}[\cite{JacZw:12}]\label{def_BCS}
Let $X$, $U$ and $Y$ be complex Hilbert spaces. We call 
\begin{equation}\label{BCS}
    \begin{aligned}
        \dot{x}(t) &= \A x(t), \quad x(0)=x_0, \\
        u(t) &= \BB x(t), \\
        y(t) &= \CC x(t),
    \end{aligned}
\end{equation}
a \emph{boundary control and observation system} on $(X,U,Y)$ if the following hold: 
\begin{enumerate}
    \item $\A : D(\A)\subset X \to X$, $\BB : D(\BB) \subset X \to U$, $\CC : D(\A) \subset X \to Y$ are linear with $D(\A) \subset D(\BB)$.
    \item The operator $A :D(A) \to X$ with $D(A) = D(\A) \cap \ker(\BB)$ and $$ Ax = \A x \quad \text{for } x \in D(A) $$
    is the generator of a $C_0$-semigroup $(T(t))_{t\geq0}$ on $X$.
   \item There exists an operator $B \in L(U, X)$ such that for all $u \in U$ we have $Bu \in D(\A)$, ${\A B \in L(U, X)}$ and $$ \BB Bu = u, \quad u \in U. $$
  \item The operator $\CC$ is bounded from the domain of $\A$ to $Y$. Here, $D(A)$ is equipped with the graph norm.
\end{enumerate}
\end{defn}

\begin{thm}[\cite{JacZw:12}]
The transfer function of the boundary control and observation system \eqref{BCS} is given by
\begin{equation*}
    \G(s) = \CC(sI - A)^{-1} (\A B - sB) + \CC B, \quad s \in \rho(A).
\end{equation*}
Here $\rho(A)$ denotes the resolvent set of the operator $A$.
For $s \in \rho(A)$ and $u_0 \in U$, $\G(s)u_0$ can also be calculated as the (unique) solution of
\begin{align*}
    s x_0 &= \A x_0, \\
    u_0 &= \BB x_0, \\
    \G(s)u_0 &= \CC x_0,
\end{align*}
with $x_0 \in D(\A)$.
\end{thm}

\begin{defn}
    The function $x : [0,\tau] \to X$ is a \emph{classical solution of the boundary control and observation system} \eqref{BCS} on $[0,\tau]$ if $x$ is a continuously differentiable function, $x(t) \in D(A)$ for all $t \in [0,\tau]$, and $x(t)$ satisfies \eqref{BCS} for all $t \in [0,\tau]$, (\cite{CurtZwa:20}).
\end{defn}

\begin{defn}
    The boundary control and observation system \eqref{BCS} on $(X,U,U)$ is \emph{impedance passive} if, for all $t>0$ and for every $x_0 \in X$ and $u \in C^2(\R^+;U)$ with ${x_0 - u(0) \in D(A)}$ the solution $(x,y)$ of \eqref{BCS} satisfies
    $$ \frac{d}{dt} \Vert x(t)\Vert^2 \leq 2\Re \langle y(t),u(t) \rangle. $$
\end{defn}

Next we show that system \eqref{eqn:system1} is an impedance passive boundary control and observation system. 
Therefore, we assume that the matrices $P_0$, $P_2$, $\HH$, $\begin{bsmallmatrix}
    W_{B,1} \\ W_{B,2}
\end{bsmallmatrix}$ and $W_C$ satisfy the assumptions made in the introduction. We define the spaces $X$, $U$, $Y$ by $X = L^2([a,b];\mathbb F^n)$ and $U=Y =\F^m$. Let $\tau : H^2([a,b];\F^n) \to \F^{4n}$, the trace operator, be given by
\begin{equation*}
    \tau(x) =
    \begin{pmatrix}
        x(b) & x'(b) &
        x(a) & x'(a)
    \end{pmatrix}^\top.
\end{equation*}
We define
\begin{equation*}
    \begin{aligned} 
         \A x &=  \left(P_2 \frac{\partial^2}{\partial \xi^2} +P_0\right)\HH x, \\
        \BB x &=  W_{B,1} \HH \tau (x) ,\\
         \CC x &=  W_C \HH \tau (x) ,\\
        D(\A) &= D(\BB)=\left\{ x \in H^2([0,1];\F^n) \,\vert\, W_{B,2}\HH\tau(x) =0 \right\}.
    \end{aligned}   
\end{equation*} 

\begin{prop}[\cite{Aug:16}]
    The port-Hamiltonian system \eqref{eqn:system1} is an impedance passive boundary control and observation system.
\end{prop}

\begin{defn}
    The boundary control and observation system $(\A,\BB,\CC)$ on $(U,X,Y)$ is called \emph{well-posed}, if for some $t > 0$, there exists $\kappa_t>0$ such that for every $x_0 \in D(A)$ and $u \in C^2(\R^+;U)$ with ${u(0)=\BB x_0}$ the solution for \eqref{BCS} on $[0,t]$ fulfills
    \begin{equation*}
    \Vert x(t)\Vert^2 +\Vert y\Vert^2_{L^2((0,t);Y)} \leq \kappa_t \left( \Vert x_0\Vert^2 +\Vert u \Vert^2_{L^2((0,t);U)} \right).
    \end{equation*}
\end{defn}

The following two theorems will be useful for the proof of our main result.

\begin{thm}\label{thm:well-posedness}
    An impedance passive boundary control and observation system is well-posed if and only if its transfer function $\G$ is bounded on some vertical line in $\C^+:=\{s \in \C ~\vert\,  \Re ~s > 0 \}$, i.e.~$\G : \C^+ \longrightarrow L(U,Y)$ satisfies:
    $$ \text{There exists } r >0 \text{ such that } \sup_{\omega \in \R} \| \G(r + i\omega)\| <\infty. $$
\end{thm}
\textbf{Proof:}
    Since every boundary control and observation system is a system node (\cite{Villegas:07}), we can apply Staffans (\cite{Sta:02}) which states that the well-posedness of an impedance passive system node is equivalent to the fact that the transfer function is bounded on some vertical line in the right half plane. 
\qed

\begin{thm}[\cite{Wei:94,Sta:05}]\label{thm:feedback}
Let $\G$ be the transfer function of the well-posed boundary control and observation system
     \begin{align*}
         \dot{x}(t) &= \A x(t) , \\
         u(t) &= \mathcal{B} x(t) , \\
         y(t) &= \mathcal{C}x(t),
    \end{align*}
Let $F$ be a bounded linear operator from $Y$ to $U$ such that the inverse of $I+\G(s)F$ $\left(\text{or } I+F\G(s)\right)$ exists and is bounded in $s$ in some right half-plane. Then the closed-loop system is again well-posed, that is, the boundary control and observation system
    \begin{align*}
         \dot{x}(t) &= \A x(t) , \\
         v(t) &=(\mathcal{B} + F\mathcal{C})x(t) , \\
         y(t) &= \mathcal{C}x(t),
    \end{align*}
    is again well-posed.
    The operator $F$ is called an \emph{admissible feedback operator}. 
\end{thm}

\section{Well-posedness of the port-Hamiltonian system}

Throughout this section we assume that the matrices $P_0$, $P_2$, $\HH$, $\begin{bsmallmatrix}
    W_{B,1} \\ W_{B,2}
\end{bsmallmatrix}$ and $W_C$ satisfy the assumptions made in the introduction. Further, we assume $m=2n$, implying that only $W_{B,1}$ exists and $W_{B,2}$ does not. We define
\begin{align}\label{eqn:usys}
    u_s(t) = \begin{bmatrix} 
           i \tfrac{\partial}{\partial \xi} \HH x(b,t), \\
           i \tfrac{\partial}{\partial \xi} \HH x(a,t) 
        \end{bmatrix} , 
    \quad 
    y_s(t) = \begin{bmatrix} 
                \HH x(b,t) \\ 
                -\HH x(a,t) 
            \end{bmatrix}. 
\end{align}
Both the input $u$ and the output $y$ can be expressed as linear combinations (or permutations) of the components of $u_s$ and $y_s$. Therefore, there exist uniquely determined matrices $B_1$, $B_2$, $C_1$, $C_2\in \F^{2n\times 2n}$ such that the port-Hamiltonian system \eqref{eqn:system1} can be equivalently written as 
\begin{equation}\label{eqn:system2}
\begin{aligned}
     \frac{\partial x}{\partial t}(\xi,t) &= \left( P_2 \frac{\partial^2}{\partial \xi^2} + P_0 \right)\HH x(\xi,t), \, \xi \in [a,b] \\
        u(t) &= B_1 u_s(t) +  B_2
        y_s(t),\\
        y(t) &= C_1  u_s(t) + C_2  y_s(t).
    \end{aligned}
\end{equation}
Next, we formulate the main result of this paper.
  
\begin{thm}\label{main_result}
The port-Hamiltonian system \eqref{eqn:system1} is well-posed if and only if the matrix $B_1$ is invertible.
\end{thm}
\begin{rem}\label{remark}
If $0 < m \leq 2n$, that is, when we are dealing with quantities that are set to zero, such as with the boundary condition $W_{B,2} \HH \tau(x) =0$, the sufficient condition still holds. This follows from the fact that we can always introduce an alternative input $v(t)=W_{B,2} \HH \tau(x)$ and prove that the system with the extended input $\begin{bsmallmatrix}
        u \\ v
    \end{bsmallmatrix}$
is well-posed. We can then choose any input that is (locally) square integrable. In particular, $v \equiv 0$.
\end{rem}

To prove Theorem \ref{main_result}, we will need the following lemmas.
\begin{lem}\label{Lemma_sin}
    For $r \in (0,\infty)$ and $x \in \R \backslash \{0\}$, we define ${\alpha(x)=\frac{r}{x}+ix}$. Then 
    \begin{equation*}
        \vert \alpha(x) \sinh(\alpha(x))\vert \geq r.
    \end{equation*} 
\end{lem}
\textbf{Proof:}
We calculate
    \begin{align*}
        \lefteqn{\vert \alpha(x) \sinh(\alpha(x))\vert^2}\\
        &= \left( \frac{r^2}{x^2} +x^2\right) \vert \sinh(\alpha(x))\vert^2 \\
        &= \left( \frac{r^2}{x^2} +x^2\right) \frac{1}{4} \vert e^{\alpha(x)} -e^{-\alpha(x)}\vert^2 \\
        &= \left( \frac{r^2}{x^2} +x^2\right) \frac{1}{2} \left( \cosh(\frac{2r}{x}) -1 + 2\sin^2(x) \right) \\
        &\geq \left( \frac{r^2}{x^2} +x^2\right) \frac{1}{2} \left( \cosh(\frac{2r}{x}) -1 \right) \\
        &\geq \left( \frac{r^2}{x^2} +x^2\right) \frac{1}{2} \frac{\left(\frac{2r}{x}\right)^2}{2} = \left( \frac{r^4}{x^4}+ r^2\right) \\
        &\geq r^2,
    \end{align*}
which concludes the proof.
\qed

\begin{lem}\label{lemma:scalar case}
Let $\mu\in \R\backslash\{0\}$. 
Then the boundary control and observation system
    \begin{align*}
        \frac{\partial x}{\partial t}(\xi,t) &= i \mu \frac{\partial^2}{\partial \xi^2} x(\xi,t), \quad \xi \in [a,b], \quad t \geq 0, \\ 
        u(t) &= \begin{bmatrix} i\mu \tfrac{\partial}{\partial \xi} x(b,t) \\[0.8ex] i\mu \tfrac{\partial}{\partial \xi} x(a,t) \end{bmatrix}, \\ 
        y(t) &= \begin{bmatrix} \mu x(b,t) \\ -\mu x(a,t) \end{bmatrix},   
    \end{align*}
with state space $X=L^2([a,b],\F)$, $U=Y=\F^2$ and Hamiltonian $ H(x)=\tfrac{1}{2}\int_a^b |\mu| |x(\xi)|^2 d\xi$ is a well-posed port-Hamiltonian system. In particular, the system is impedance passive. Moreover, the transfer function $\G_\mu$ of the system satisfies
    $$ \lim_{\Re\, s \to \infty} \G_\mu(s) =0 . $$
\end{lem}

\begin{pf}
Let $n=1$, $m=2$, $P_2=i\operatorname{sgn}(\mu)$, $P_0=0$, $\HH= |\mu|$, $W_{B,1}= \begin{bsmallmatrix}
    0 & i &0&0\\ 0&0&0&i
\end{bsmallmatrix}$ and $W_{C}= \begin{bsmallmatrix}
    1 & 0 &0&0\\ 0&0& -1 &0
\end{bsmallmatrix}$. Thus $R = \begin{bsmallmatrix}
    0 & i \\ -i&0
\end{bsmallmatrix}$, $\tilde W_B=\tfrac{1}{\sqrt{2}}\begin{bsmallmatrix}
    1 & 0 &0&i\\ -1&0&0&i
\end{bsmallmatrix}$ and $\tilde W_C=\tfrac{1}{\sqrt{2}}\begin{bsmallmatrix}
    0 & i &1&0\\ 0&i&-1&0
\end{bsmallmatrix}$. This implies 
$$
    \begin{bmatrix}  \tilde W_B\Sigma \tilde W_B^\ast &  \tilde W_B\Sigma \tilde W_C^\ast \\  \tilde W_C\Sigma \tilde W_B^\ast &  \tilde W_C\Sigma \tilde W_C^\ast \end{bmatrix}^{-1}= \begin{bmatrix}  0 & I \\ I & 0 \end{bmatrix}.
$$
Thus the system is port-Hamiltonian system and in particular impedance passive.
By Theorem \ref{thm:well-posedness}, in order to show well-posedness it remains to show that the corresponding transfer function on some vertical line in the open right half plane. For this further, it is sufficient to show that there exists a constant $M>0$ such that 
$$ \| \G_\mu(s) \| \leq \frac{M}{\Re\, s}. $$

Let $s\in \mathbb C$ with $\Re\, s>0$ and $u_0\in \F^2$.
The vector $\G_\mu(s)u_0$ is uniquely determined by (see \cite{JacZw:12})
    \begin{align*}
        s x_0(\xi) &= i\mu \frac{\partial^2}{\partial \xi^2}x_0(\xi), \quad \xi \in [a,b], \\
        u_{0} &= \begin{bmatrix} i\mu \tfrac{\partial}{\partial \xi}x_0(b) \\[0.8ex] i\mu \tfrac{\partial}{\partial \xi}x_0(a) \end{bmatrix}, \\
        \G_\mu(s)u_{0} &= \begin{bmatrix} \mu x_0(b) \\ -\mu x_0(a) \end{bmatrix}.
    \end{align*}
    The solution of $s x_0(\xi) = i\mu \frac{\partial^2}{\partial \xi^2}x_0(\xi)$ is given by 
    $$ x_0(\xi) = \alpha e^{\gamma\xi} + \beta e^{-\gamma\xi}, \quad \alpha,\beta \in \C , $$
    where $\gamma^2=-is\mu^{-1}$.
    
    Considering the boundary conditions, we find that 
    \begin{equation*}
        \begin{bmatrix}
            \alpha \\ \beta
        \end{bmatrix} = \frac{i}{\mu \gamma (e^{\gamma(b-a)} - e^{-\gamma(b-a)})}
        \begin{bmatrix}
            - e^{-\gamma a} &  e^{-\gamma b} \\
            - e^{\gamma a} &  e^{\gamma b}
        \end{bmatrix} u_0.
    \end{equation*}
    Hence, the transfer function is given by 
    \begin{equation*}
        \G_\mu(s) = \frac{i}{\gamma} 
        \begin{bmatrix}
            -\coth(\gamma(b-a)) & \sinh(\gamma(b-a))^{-1} \\
            \sinh(\gamma(b-a))^{-1} & -\coth(\gamma(b-a)) 
        \end{bmatrix}.
    \end{equation*} 
    
    Let $r \in (0,\infty)$ and $\omega \in \R$. For $s=r+i\omega$, we find $\Re \, \gamma = \frac{-r}{2 \mu \Im\,\gamma}$. Since $\G_\mu$ is a function of $s$, it does not depend on the choice of the square root of $\gamma^2$. Therefore, we can always choose a square root $\gamma$ in the right half-plane. Without loss of generality, we assume that $\Re\, \gamma > 0$ and it follows from Lemma \ref{Lemma_sin} that
    \begin{equation*}
        \vert \gamma^{-1}\sinh(\gamma(b-a))^{-1} \vert \leq \frac{2 \vert (b-a) \mu \vert}{r} . 
    \end{equation*}
    Furthermore, it holds that 
    \begin{equation*}
        \coth(\gamma(b-a))^2 = 1 + \frac{1}{\sinh(\gamma(b-a))^2} .
    \end{equation*} 
    Thus, the norm of the matrix $\G_\mu(s)$ converges to 0 as $\Re\, s \to \infty$. 
    \qed
\end{pf}

Having established these two lemmas, we can now proceed with the proof of Theorem \ref{main_result}. 

\textbf{Proof of Theorem \ref{main_result}.}
As $P_0 \HH$, seen as a multiplication operator on $X$, is a bounded operator on the state space $X$ and the well-posedness remains under bounded perturbations \cite[Lemma 13.1.14]{JacZw:12}, we may assume without loss of generality that $P_0=0$.
The proof is structured in several steps.

\textit{Step 1.} There exists an invertible matrix $Q$ and a diagonal matrix $\Delta$ such that  
\begin{equation*}
    P_2 \HH = Q^{-1} \Delta Q,
\end{equation*} 
where 
\begin{equation}\label{eqn:delta}
    \Delta = i
    \begin{bmatrix}
        \Delta^+ & 0 \\
        0 & \Delta^-
    \end{bmatrix},
\end{equation}
$\Delta^+$ contains the positive and $\Delta^-$ the negative
eigenvalues of $P_2 \HH$ on the respective diagonals.

\textit{Proof:} The statement follows from the fact that $\HH^{1/2}P_2\HH^{1/2}$ is skew-adjoint, $P_2 \HH = \HH^{-1/2}(\HH^{1/2}P_2\HH^{1/2}) \HH^{1/2}$ and $P_2 \HH$ is invertible.\qed

\textit{Step 2.} Without loss of generality we may assume that the port-Hamiltonian system \eqref{eqn:system2} is given by 
\begin{equation}\label{eqn:system3}
    \begin{aligned}
        \frac{\partial}{\partial t}
        \begin{bmatrix}
            x^+ (\xi,t) \\ 
            x^- (\xi,t)
        \end{bmatrix} &= i 
        \begin{bmatrix}
            \Delta^+ & 0 \\
            0 & \Delta^-
        \end{bmatrix} \frac{\partial^2}{\partial \xi^2} 
        \begin{bmatrix}
            x^+ (\xi,t) \\ 
            x^- (\xi,t)
        \end{bmatrix},  \\  
        u(t) &=    \tilde B_1 u_s(t) +    \tilde B_2 y_s(t),\\
        y(t) &= \tilde C_1  u_s(t) + \tilde C_2  y_s(t),\\
        u_s(t) &= 
        \begin{bmatrix} 
            i \Delta^+ \tfrac{\partial}{\partial \xi} x^+(b,t) \\
            i\Delta^- \tfrac{\partial}{\partial \xi} x^-(b,t) \\
            i \Delta^+ \tfrac{\partial}{\partial \xi} x^+(a,t) \\
            i \Delta^- \tfrac{\partial}{\partial \xi} x^-(a,t) 
        \end{bmatrix}, \\
        y_s(t) &= 
        \begin{bmatrix} 
             \Delta^+ x^+(b,t) \\ 
             \Delta^- x^-(b,t) \\
             -\Delta^+ x^+(a,t) \\
             -\Delta^- x^-(a,t) 
        \end{bmatrix},
         \end{aligned}
    \end{equation}
where $\Delta^+$ contains the positive and $\Delta^-$ the negative
eigenvalues of $P_2 \HH$ on the respective diagonals. Moreover, $\tilde B_1= B_1 P_2^{-1}Q^{-1}$. Thus $B_1$ is invertible if and only if $\tilde B_1$ is invertible.

\textit{Proof:} Applying the state space transformation
$$\left[\begin{smallmatrix}
            x^+  \\ 
            x^- 
        \end{smallmatrix}\right]=Q x$$
to the port-Hamiltonian system \eqref{eqn:system2} with $P_0=0$ we obtain
the port-Hamiltonian system \eqref{eqn:system3}. After the bounded state transformation, the matrices $B_1$, $B_2$, $C_1$, and $C_2$ are transformed but remain the same.
\qed

\textit{Step 3:} The port-Hamiltonian system \begin{equation}\label{eqn:phs_diagonal}
    \begin{aligned}
        \frac{\partial}{\partial t}
        \begin{bmatrix}
            x^+ (\xi,t) \\ 
            x^- (\xi,t)
        \end{bmatrix} &= i 
        \begin{bmatrix}
            \Delta^+ & 0 \\
            0 & \Delta^-
        \end{bmatrix} \frac{\partial^2}{\partial \xi^2} 
        \begin{bmatrix}
            x^+ (\xi,t) \\ 
            x^- (\xi,t)
        \end{bmatrix},  \\ 
                 u_s(t) &= 
        \begin{bmatrix} 
            i \Delta^+ \tfrac{\partial}{\partial \xi} x^+(b,t) \\
            i\Delta^- \tfrac{\partial}{\partial \xi} x^-(b,t) \\
            i \Delta^+ \tfrac{\partial}{\partial \xi} x^+(a,t) \\
            i \Delta^- \tfrac{\partial}{\partial \xi} x^-(a,t) 
        \end{bmatrix}, \\
        y_s(t) &= 
        \begin{bmatrix} 
              \Delta^+ x^+(b,t) \\ 
              \Delta^- x^-(b,t) \\
              -\Delta^+ x^+(a,t) \\
              -\Delta^- x^-(a,t) 
        \end{bmatrix},
         \end{aligned}
    \end{equation}
with input $u_s$ and output $y_s$ is well-posed and the corresponding transfer function $\G_s$ satisfies
\begin{equation*}
    \lim_{\Re\, s \to \infty} \G_s(s) = 0 .
\end{equation*}

\textit{Proof:}
As the port-Hamiltonian system \eqref{eqn:system3}
consists of copies of systems studied in Lemma \ref{lemma:scalar case}, the statement follows directly.
\hfill\qed

\textit{Step 4.} 
If $B_1$ is invertible, or equivalently, if $\tilde B_1$ is invertible, then the port-Hamiltonian system \eqref{eqn:system3} is well-posed.
 
\begin{figure}[ht]
    \begin{center}
    \begin{tikzpicture}[scale=1]
        \node[draw] (I) at (0.5,0) {$\ \G_s(s)\ $};
        \node[draw] (II) at (2.8,0) {$ \ \tilde C_2 \ $};
        \node[draw] (III) at (-2.3,0) {$  \tilde B_1
       ^{-1} \ $};
        \node[draw] (ID) at (1,1) {$ \ \tilde C_1 \ $}; 
        \node[draw] (IV) at (0.6,-1) {$ \ \tilde B_1
       ^{-1}\tilde B_2 \ $};
        \fill (-1,0) circle (0.9pt); 
        \fill[black] (3.4,1) circle(0.9pt); 
        \fill[black] (2,0) circle(0.9pt); 
        \draw (3.8,1) node[above]{$y$}; \draw[->,>=latex] (3.3,1)--(4.2,1); 
        \draw (1.5,0) node[above]{$y_s$};
        \draw (-3.2,0) node[above]{$u$}; \draw[->,>=latex] (-3.6,0)--(-2.8,0); 
        \draw (-0.6,0) node[above]{$u_s$};  
        \draw[->,>=latex] (I)--(II); \draw[->,>=latex] (III)--(I); \draw (ID)-|(-1,0); \draw (ID)-|(3.4,0); \draw (II)--(3.4,0); \draw (IV)-|(2,0);
        \draw[->,>=latex] (IV)-|(-1,0);
    \end{tikzpicture}
    \caption{The port-Hamiltonian system  \eqref{eqn:system3}.}
    \label{figg}
    \end{center}    
    \end{figure}
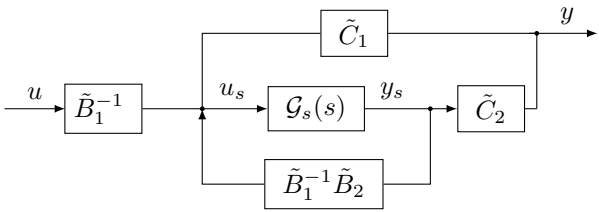

 \textit{Proof:} In this case we can write $u_s$ as  
    \begin{equation*}
        u_s(t) =  \tilde B_1^{-1} u(t) - \tilde B_1^{-1} \tilde B_2 y_s(t).
    \end{equation*} 
   Thus, the system contains a feedback loop with the feedback defined by $-\tilde B_1^{-1}\tilde B_2$ (see Figure \ref{figg}). 
    Since $\displaystyle\lim_{\Re\, s \to \infty} \G_s(s)=0$, the matrix $\tilde B_1^{-1} \tilde B_2$ defines an admissible feedback operator and thus by Theorem \ref{thm:feedback} the closed loop system is a well-posed boundary system. 

    As the closed-loop system is the port-Hamiltonian system \eqref{eqn:system3}, the statement follows.\qed 

 \textit{Step 5:} If the port-Hamiltonian system \eqref{eqn:system3} is well-posed, then the matrix $B_1$ is invertible. 

 \textit{Proof:} Assume that $B_1$ is not invertible. We denote by $\G$ the transfer function of the port-Hamiltonian system \eqref{eqn:system3}. 
  Using
    \begin{align*}
        u(t) &= \tilde B_1 u_s(t) + \tilde B_2 y_s(t), \\
        y(t) &= \tilde C_1 u_s(t) + \tilde C_2 y_s(t),
    \end{align*}
    it follows that
    \begin{equation}\label{eqn:transfer_fct}
        \G(s) \left( \tilde B_1 + \tilde B_2 \G_s(s) \right) u_{s,0} = \left( \tilde C_1 + \tilde C_2 \G_s(s) \right) u_{s,0}
    \end{equation}
    for $\Re \,s>\alpha$.

    We choose $u_{s,0} \in \ker \tilde B_1$. This choice implies that $\tilde C_1 u_{s,0} \neq 0$, as the assumption $\begin{bsmallmatrix} W_{B,1} \\ W_C \end{bsmallmatrix}$ invertible implies that the matrix $\begin{bsmallmatrix}  B_1 &  B_2 \\  C_1 & C_2 \end{bsmallmatrix}$ is invertible as well, or equivalently, $\begin{bsmallmatrix} \tilde B_1 & \tilde B_2 \\ \tilde C_1 & \tilde C_2 \end{bsmallmatrix}$ is invertible. However, taking the limit as $\Re\, s \to\infty$ in equation \eqref{eqn:transfer_fct}, we observe that in this limit $\tilde C_1 u_{s,0}=0$, which is a contradiction to $\tilde C_1 u_{s,0} \not= 0$.
    \qed

\section{Examples}

\begin{exmp}[\cite{Aug:16}]
We study the Schrö\-dinger equation of a free particle of mass $m > 0$ moving in a one-dimensional direction. The Schrö\-dinger equation is then given by
\begin{equation*} 
    i \hbar \frac{\partial}{\partial t}\psi(\xi,t) = - \frac{\hbar^2}{2m} \frac{\partial^2}{\partial \xi^2} \psi(\xi,t), \quad \xi \in [0,1], \ t\geq 0 ,
\end{equation*}  
where $\psi(\xi,t)$ denotes the wave function and $\hbar >0$ is the reduced Planck constant. We choose the following boundary control and observation
$$ 
    u(t) = \begin{bmatrix}
        \frac{\partial}{\partial \xi}\psi(1,t) \\
        i \frac{\hbar}{2m} \frac{\partial}{\partial \xi}\psi(0,t)
    \end{bmatrix},
    \quad
    y(t) = \begin{bmatrix}
        \psi(0,t) \\ \psi(1,t)
    \end{bmatrix} .  $$
By letting
    $$ P_2=i, \quad P_0=0, \text{ and } \HH= \frac{\hbar}{2m}, $$
we obtain the port-Hamiltonian formulation for the Schrö\-dinger equation, and for this choice, the matrix $B_1$ is given by 
$$ B_1 =  \begin{bmatrix} 1 & 0 \\ 0 & i \end{bmatrix}, $$
which is invertible. Thus, the system is well-posed.\qed
\end{exmp}

\begin{exmp}[\cite{HumLassi:18}]\label{example}
    The Euler-Bernoulli beam equation is given on the spatial interval $[0,1]$ by
    $$ \rho \frac{\partial^2 w}{\partial t^2}(\xi,t) + \frac{\partial^2}{\partial \xi^2}\left( EI\frac{\partial^2 w}{\partial \xi^2}(\xi,t)\right)=0. $$
    We have already formally expressed the Euler-Bernoulli beam equation in the form port-Hamiltonian system \eqref{eqn:system1}, and we define the control and observation operators as follows 
    \begin{equation*}
        u(t) =  
        \begin{bmatrix}
            \frac{\partial}{\partial \xi}\HH_1 x_1(0,t) \\ \HH_1 x_1(0,t) \\ \HH_2 x_2(1,t) \\ \frac{\partial}{\partial \xi}\HH_2 x_2(1,t)   
        \end{bmatrix},
        \quad
        y(t) =  
        \begin{bmatrix}
            -\HH_2 x_2(0,t) \\ \frac{\partial}{\partial \xi}\HH_2 x_2(0,t) \\ 
            -\HH_1 x_1(1,t) \\ \frac{\partial}{\partial \xi}\HH_1 x_1(1,t)  
        \end{bmatrix},
    \end{equation*} 
   where $\HH_1 = \frac{1}{\rho}$ and $\HH_2=EI$.

   For this case, we have  
  $$     B_1 = -i \begin{bmatrix}
       0 & 0 & 1 & 0 \\
       0 & 0 & 0 & 0 \\
       0 & 0 & 0 & 0 \\
       0 & 1 & 0 & 0
   \end{bmatrix}. $$
   
   Since the matrix $B_1$ is not invertible, the Euler-Bernoulli beam equation with this choice of boundary control and observation is not well-posed. \qed
\end{exmp}

\begin{exmp}[\cite{GuoWang:05}]
    We consider a beam with a roller support at one end $\xi=0$ and a free end at $\xi =1$, at which a control force applied based on the angular velocity. Hence, the boundary conditions are defined as follows 
    \begin{align*}
    \frac{\partial w}{\partial \xi}(0,t)=0 &, \quad EI \frac{\partial^3 w}{\partial \xi^3}(0,t)=0 ,  \\
    EI \frac{\partial^3 w}{\partial \xi^3}(1,t)=0 &, \quad \frac{\partial^2 w}{\partial t \partial \xi}(1,t) =u(t). 
    \end{align*}
    The boundary control are given by 
    \begin{equation*}
        \begin{bmatrix}
            0 \\ I
        \end{bmatrix}
        u(t) =  
        \begin{bmatrix}
            \tfrac{\partial}{\partial \xi}\HH_1 x_1(0,t) \\ 
            \tfrac{\partial}{\partial \xi}\HH_2 x_2 (0,t) \\ 
            \tfrac{\partial}{\partial \xi}\HH_2 x_2(1,t)   \\
            \tfrac{\partial}{\partial \xi}\HH_1 x_1(1,t) 
        \end{bmatrix},
    \end{equation*}
    and the boundary observation
    $$ y(t) =  - \HH_2 x_2(1,t), $$

   where $\HH_1 = \frac{1}{\rho}$ and $\HH_2=EI$. Here, we have  
   \begin{equation*}
    B_1 = -i \begin{bmatrix}
       0 & 0 & 1 & 0 \\
       0 & 0 & 0 & 1 \\
       0 & 1 & 0 & 0 \\
       1 & 0 & 0 & 0
      \end{bmatrix}.
   \end{equation*}
   
   The matrix $B_1$ is invertible, we conclude by Remark \ref{remark} that with this choice of boundary conditions the Euler-Bernoulli beam is well-posed. 
   \qed
\end{exmp}

\section{Conclusion}
We considered a class of port-Hamiltonian systems with boundary control and observation, and provided a necessary and sufficient condition to ensure their well-posedness. This condition can be verified through a simple matrix check. Initially, we start with the system \eqref{eqn:system2} without assuming it is well-posed in the sense of a $C_0$-semigroup, but we derive this property from the system. We find that if $B_1$ is singular, at least one of the four properties of the well-posedness is lost, potentially including the semigroup property.
As a main application, we analyzed the one-dimensional Euler-Bernoulli beam with boundary controls and observations, proving that well-posedness may not hold even when the physical parameters are constant. Future work will focus on extending this result to more general port-Hamiltonian systems, particularly in cases where the Hamiltonian density $\HH$ depends on the spatial variable.

\bibliography{ifacconf}             

\end{document}